\documentclass[12pt]{article}

\usepackage{amsmath,amssymb,amsthm,amscd,a4wide}

\begin{document}

\newcommand{\End}{{\rm{End}\ts}}
\newcommand{\Hom}{{\rm{Hom}}}
\newcommand{\Mat}{{\rm{Mat}}}
\newcommand{\ch}{{\rm{ch}\ts}}
\newcommand{\chara}{{\rm{char}\ts}}
\newcommand{\diag}{ {\rm diag}}
\newcommand{\non}{\nonumber}
\newcommand{\wt}{\widetilde}
\newcommand{\wh}{\widehat}
\newcommand{\ot}{\otimes}
\newcommand{\la}{\lambda}
\newcommand{\La}{\Lambda}
\newcommand{\De}{\Delta}
\newcommand{\al}{\alpha}
\newcommand{\be}{\beta}
\newcommand{\ga}{\gamma}
\newcommand{\Ga}{\Gamma}
\newcommand{\ep}{\epsilon}
\newcommand{\ka}{\kappa}
\newcommand{\vk}{\varkappa}
\newcommand{\si}{\sigma}
\newcommand{\vp}{\varphi}
\newcommand{\de}{\delta}
\newcommand{\ze}{\zeta}
\newcommand{\om}{\omega}
\newcommand{\ee}{\epsilon^{}}
\newcommand{\su}{s^{}}
\newcommand{\hra}{\hookrightarrow}
\newcommand{\ve}{\varepsilon}
\newcommand{\vs}{\varsigma}
\newcommand{\ts}{\,}
\newcommand{\vac}{\mathbf{1}}
\newcommand{\di}{\partial}
\newcommand{\qin}{q^{-1}}
\newcommand{\tss}{\hspace{1pt}}
\newcommand{\Sr}{ {\rm S}}
\newcommand{\U}{ {\rm U}}
\newcommand{\BL}{ {\overline L}}
\newcommand{\BE}{ {\overline E}}
\newcommand{\BP}{ {\overline P}}
\newcommand{\AAb}{\mathbb{A}\tss}
\newcommand{\CC}{\mathbb{C}\tss}
\newcommand{\KK}{\mathbb{K}\tss}
\newcommand{\QQ}{\mathbb{Q}\tss}
\newcommand{\SSb}{\mathbb{S}\tss}
\newcommand{\ZZ}{\mathbb{Z}\tss}
\newcommand{\X}{ {\rm X}}
\newcommand{\Y}{ {\rm Y}}
\newcommand{\Z}{{\rm Z}}
\newcommand{\Ac}{\mathcal{A}}
\newcommand{\Dc}{\mathcal{D}}
\newcommand{\Lc}{\mathcal{L}}
\newcommand{\Mc}{\mathcal{M}}
\newcommand{\Pc}{\mathcal{P}}
\newcommand{\Qc}{\mathcal{Q}}
\newcommand{\Tc}{\mathcal{T}}
\newcommand{\Sc}{\mathcal{S}}
\newcommand{\Bc}{\mathcal{B}}
\newcommand{\Ec}{\mathcal{E}}
\newcommand{\Fc}{\mathcal{F}}
\newcommand{\Hc}{\mathcal{H}}
\newcommand{\Uc}{\mathcal{U}}
\newcommand{\Vc}{\mathcal{V}}
\newcommand{\Wc}{\mathcal{W}}
\newcommand{\Yc}{\mathcal{Y}}
\newcommand{\Ar}{{\rm A}}
\newcommand{\Br}{{\rm B}}
\newcommand{\Ir}{{\rm I}}
\newcommand{\Fr}{{\rm F}}
\newcommand{\Hr}{{\rm H}}
\newcommand{\Jr}{{\rm J}}
\newcommand{\Or}{{\rm O}}
\newcommand{\GL}{{\rm GL}}
\newcommand{\Spr}{{\rm Sp}}
\newcommand{\Rr}{{\rm R}}
\newcommand{\Zr}{{\rm Z}}
\newcommand{\gl}{\mathfrak{gl}}
\newcommand{\middd}{{\rm mid}}
\newcommand{\ev}{{\rm ev}}
\newcommand{\Pf}{{\rm Pf}}
\newcommand{\Norm}{{\rm Norm\tss}}
\newcommand{\oa}{\mathfrak{o}}
\newcommand{\spa}{\mathfrak{sp}}
\newcommand{\osp}{\mathfrak{osp}}
\newcommand{\g}{\mathfrak{g}}
\newcommand{\h}{\mathfrak h}
\newcommand{\n}{\mathfrak n}
\newcommand{\z}{\mathfrak{z}}
\newcommand{\Zgot}{\mathfrak{Z}}
\newcommand{\p}{\mathfrak{p}}
\newcommand{\sll}{\mathfrak{sl}}
\newcommand{\agot}{\mathfrak{a}}
\newcommand{\qdet}{ {\rm qdet}\ts}
\newcommand{\Ber}{ {\rm Ber}\ts}
\newcommand{\HC}{ {\mathcal HC}}
\newcommand{\cdet}{ {\rm cdet}}
\newcommand{\tr}{ {\rm tr}}
\newcommand{\gr}{ {\rm gr}}
\newcommand{\str}{ {\rm str}}
\newcommand{\loc}{{\rm loc}}
\newcommand{\Gr}{{\rm G}}
\newcommand{\sgn}{ {\rm sgn}\ts}
\newcommand{\ba}{\bar{a}}
\newcommand{\bb}{\bar{b}}
\newcommand{\bi}{\bar{\imath}}
\newcommand{\bj}{\bar{\jmath}}
\newcommand{\bk}{\bar{k}}
\newcommand{\bl}{\bar{l}}
\newcommand{\hb}{\mathbf{h}}
\newcommand{\Sym}{\mathfrak S}
\newcommand{\fand}{\quad\text{and}\quad}
\newcommand{\Fand}{\qquad\text{and}\qquad}
\newcommand{\For}{\qquad\text{or}\qquad}
\newcommand{\OR}{\qquad\text{or}\qquad}

\renewcommand{\theequation}{\arabic{section}.\arabic{equation}}

\newtheorem{thm}{Theorem}[section]
\newtheorem{lem}[thm]{Lemma}
\newtheorem{prop}[thm]{Proposition}
\newtheorem{cor}[thm]{Corollary}
\newtheorem{conj}[thm]{Conjecture}
\newtheorem*{mthm}{Main Theorem}
\newtheorem*{mthma}{Theorem A}
\newtheorem*{mthmb}{Theorem B}
\newtheorem*{mthmc}{Theorem C}

\theoremstyle{definition}
\newtheorem{defin}[thm]{Definition}

\theoremstyle{remark}
\newtheorem{remark}[thm]{Remark}
\newtheorem{example}[thm]{Example}

\newcommand{\bth}{\begin{thm}}
\renewcommand{\eth}{\end{thm}}
\newcommand{\bpr}{\begin{prop}}
\newcommand{\epr}{\end{prop}}
\newcommand{\ble}{\begin{lem}}
\newcommand{\ele}{\end{lem}}
\newcommand{\bco}{\begin{cor}}
\newcommand{\eco}{\end{cor}}
\newcommand{\bde}{\begin{defin}}
\newcommand{\ede}{\end{defin}}
\newcommand{\bex}{\begin{example}}
\newcommand{\eex}{\end{example}}
\newcommand{\bre}{\begin{remark}}
\newcommand{\ere}{\end{remark}}
\newcommand{\bcj}{\begin{conj}}
\newcommand{\ecj}{\end{conj}}

\newcommand{\bal}{\begin{aligned}}
\newcommand{\eal}{\end{aligned}}
\newcommand{\beq}{\begin{equation}}
\newcommand{\eeq}{\end{equation}}
\newcommand{\ben}{\begin{equation*}}
\newcommand{\een}{\end{equation*}}

\newcommand{\bpf}{\begin{proof}}
\newcommand{\epf}{\end{proof}}

\def\beql#1{\begin{equation}\label{#1}}

\title{\Large\bf Invariants of the vacuum module associated with\\
the Lie superalgebra $\gl(1|1)$}

\author{{A. I. Molev\quad and\quad
E. E. Mukhin}}

\date{} 
\maketitle

\vspace{25 mm}

\begin{abstract}
We describe the algebra of invariants of the vacuum module associated
with an affinization of the Lie superalgebra $\gl(1|1)$.
We give a formula for its Hilbert--Poincar\'{e} series
in a fermionic (cancellation-free) form which turns out to coincide with
the generating function of the plane partitions over
the $(1,1)$-hook. Our arguments are based on a super version of the
Beilinson--Drinfeld--Ra\"{i}s--Tauvel theorem which we prove by producing
an explicit basis of invariants of the symmetric algebra of polynomial currents
associated with $\gl(1|1)$. We identify the invariants with affine supersymmetric
polynomials via a version of the Chevalley theorem.

%

\end{abstract}


\vspace{30 mm}

\noindent
School of Mathematics and Statistics\newline
University of Sydney,
NSW 2006, Australia\newline
alexander.molev@sydney.edu.au

\vspace{7 mm}

\noindent
Department of Mathematical Sciences\newline
Indiana University -- Purdue University Indianapolis\newline
402 North Blackford St, Indianapolis, IN 46202-3216, USA\newline
mukhin@math.iupui.edu

%

\newpage

\hfill {\sl Dedicated to Rodney Baxter on his 75th birthday}

\section{Introduction}
\label{sec:int}
\setcounter{equation}{0}

Suppose that $\g$ is a
simple Lie algebra over $\CC$ and $\ka\in\CC$.
The vacuum module $V_{\ka}(\g)$ at the level $\ka$
over the affine Kac--Moody algebra $\wh\g$
has a vertex algebra structure. The center of this vertex algebra
is trivial unless the level is
{\em critical}; this is a unique value of $\ka$ depending on the normalization of the
invariant symmetric bilinear form on $\g$.
In this case the center $\z(\wh\g)$ is an algebra of polynomials in infinitely many
variables as described by a theorem of Feigin and Frenkel~\cite{ff:ak};
see also~\cite{f:lc}.
One may regard the Feigin--Frenkel
center $\z(\wh\g)$ as a commutative subalgebra of the universal
enveloping algebra $\U\big(t^{-1}\g[t^{-1}]\big)$. This leads to
connections with the Gaudin model
and to constructions of commutative subalgebras
of $\U(\g)$ and its tensor powers; see \cite{ffr:gm}, \cite{ffr:oi}, \cite{fft:gm} and \cite{r:si}.
Explicit constructions of generators of the Feigin--Frenkel
center for the classical Lie algebras $\g$ were given in \cite{cm:ho}, \cite{ct:qs}
and \cite{m:ff}. The Harish--Chandra images of the generators
as elements of the corresponding classical $\Wc$-algebras
were calculated in \cite{mm:yc} with the use of the Yangian characters.

As explained in Kac~\cite[Sec.~4.7]{k:va}, the construction of the vertex algebra
$V_{\ka}(\g)$ can be extended to any Lie superalgebra $\g$ equipped with
an invariant supersymmetric bilinear form. In the case of general linear Lie
superalgebras $\g=\gl(m|n)$, constructions of several families of elements
of the center $\z\big(\wh\gl(m|n)\big)$ at the critical level
were given in \cite{mr:mm}. It was
conjectured there that each of the families generates the center.
The main result of the present paper is a proof
of the conjecture in the case $m=n=1$. We believe this result will be
a key ingredient for the proof of the conjecture for arbitrary $m$ and $n$.
To give its precise formulation, we will recall the construction of \cite{mr:mm}
in more detail.

\subsection{Segal--Sugawara vectors for $\gl(m|n)$}
\label{subsec:ssvec}

For $\g=\gl(m|n)$ consider the Lie superalgebra
$\wh\g=\g[t,t^{-1}]\oplus\CC K$
with the commutation relations\footnote{The element $K$ corresponds
to $K'$ in the corrected version of \cite{mr:mm} in {\tt arXiv:0911.3447v4}.}
\begin{align}
\non
\big[E_{ij}[r],E_{kl}[s\tss]\tss\big]
=\de_{kj}\ts E_{i\tss l}[r+s\tss]
{}&-\de_{i\tss l}\ts E_{kj}[r+s\tss](-1)^{(\bi+\bj)(\bk+\bl)}\\
\label{commrel}
{}&+K\Big((n-m)\tss\de_{kj}\tss\de_{i\tss l}(-1)^{\bi}
+\de_{ij}\tss\de_{kl}(-1)^{\bi+\bk}\Big)
\ts r\tss\de_{r,-s},
\end{align}
the element $K$ is even and central, and
we set $E_{ij}[r]=E_{ij}t^r$. The $\ZZ_2$-degree (or parity) of the element $E_{ij}[r]$
is $\bi+\bj$, where $\bi=0$ for $1\leqslant i\leqslant m$ and
$\bi=1$ for $m+1\leqslant i\leqslant m+n$.

The {\it vacuum module}
$V_{\text{cri}}(\g)$ at the critical level over $\wh\g$ is defined as the quotient
of the universal enveloping algebra
$\U(\wh\g)$ by the left ideal generated by
$\g[t]$ and $K-1$. It possesses a vertex algebra
structure; see e.g. \cite{k:va}.
The center of the vertex algebra $V_{\text{cri}}(\g)$
is defined by
\ben
\z(\wh\g)
=\{S\in V_{\text{cri}}(\g)\ |\ \g[t]\ts S=0\}.
\een
Elements of $\z(\wh\g)$
are called {\it Segal--Sugawara vectors\/}.
The axioms of the vertex algebra imply that the center
is a commutative associative superalgebra and it can be
identified with
a commutative subalgebra of $\U\big(t^{-1}\g[t^{-1}]\big)$.
The application of the state-field correspondence map to
a Segal--Sugawara vector yields a field whose
Fourier coefficients are {\it Sugawara operators\/} commuting
with the action of $\wh\g$. In particular, they form
a commuting family of $\wh\g$-endomorphisms
of Verma modules over $\wh\g$ at the critical level.

The main results of \cite{mr:mm} include
an explicit construction
of several families of Segal--Sugawara vectors.
To reproduce them,
we will use
the extended Lie superalgebra $\wh\g\oplus\CC\tau$,
where the element $\tau$ is even and
\ben
\big[\tau,E_{ij}[r]\tss\big]=-r\ts E_{ij}[r-1],\qquad
\big[\tau,K\big]=0.
\een
Consider the square matrix $Z=[Z_{ij}]$ with
\ben
Z_{ij}=\de_{ij}\tau+E_{ij}[-1](-1)^{\bi}
\een
with the entries in the universal enveloping algebra $\U$
for $\wh\g\oplus\CC\tau$.
We will identify the matrix $Z$ with an element
of the tensor product superalgebra $\End\CC^{m|n}\ot \U$ by
\ben
Z=\sum_{i,j=1}^{m+n}e_{ij}\ot Z_{ij}(-1)^{\bi\tss\bj+\bj},
\een
where the $e_{ij}$ denote the standard matrix units.
Taking multiple tensor products
\beql{tenpr}
\End\CC^{m|n}\ot\dots\ot\End\CC^{m|n}\ot\U
\eeq
with $k$ copies of $\End\CC^{m|n}$,
for any $a=1,\dots,k$ we will write $Z_a$
for the matrix $Z$ corresponding to the $a$-th copy of
the endomorphism superalgebra
so that the components in all remaining
copies are the identity matrices. The symmetric group $\Sym_k$
acts naturally on the tensor product
space $(\CC^{m|n})^{\ot k}$. We let $H_k$ and $A_k$ denote
the respective images of the normalized
symmetrizer and antisymmetrizer
\ben
\frac{1}{k!}\sum_{\si\in\Sym_k}\si\in\CC[\Sym_k],
\qquad\qquad
\frac{1}{k!}\sum_{\si\in\Sym_k}\sgn\si\cdot\si\in\CC[\Sym_k]
\een
in \eqref{tenpr}. Recall that the {\it supertrace\/}
of an even matrix $X=[X_{ij}]$ with
entries in a superalgebra is defined by
\ben
\str\ts X
=\sum_{i=1}^{m+n} X_{ii}(-1)^{\bi}.
\een
Furthermore, if $X$ is invertible,
denote by $X^{\tss\prime}_{ij}$ the matrix elements of its inverse
so that $X^{-1}=[X^{\tss\prime}_{ij}]$. The
(noncommutative) {\it Berezinian\/} of $X$ is defined as the product of two
determinants
\ben
\Ber X=\sum_{\si\in\Sym_m}\sgn\si\cdot
X_{\si(1)1}\dots X_{\si(m)m}
\ts
\sum_{\tau\in\Sym_n}\sgn\tau\cdot
X^{\tss\prime}_{m+1,m+\tau(1)}\dots X^{\tss\prime}_{m+n,m+\tau(n)}.
\een

By the results of \cite[Sec.~3.1]{mr:mm},
all the coefficients
$s_{k\tss l},b_{k\tss l}, h_{k\tss l}\in\U\big(t^{-1}\g[t^{-1}]\big)$
in the expansions
\begin{align}
\str\ts Z^k&={}
s_{k\tss 0}\ts\tau^{k}+s_{k1}\ts\tau^{k-1}
+\dots+s_{kk},
\non\\[0.4em]
\str^{}_{1,\dots,k}\ts A_k Z_1\dots Z_k&={}
b_{k\tss 0}\tau^k+b_{k1}\tau^{k-1}+\dots+b_{kk},
\non\\[0.4em]
\str^{}_{1,\dots,k}\ts H_k Z_1\dots Z_k&={}
h_{k\tss 0}\tau^k+h_{k1}\tau^{k-1}+\dots+h_{kk}
\label{hkl}
\end{align}
are Segal--Sugawara vectors. Moreover, these coefficients can be expressed in terms of the
Berezinian through the identities
\begin{align}\non
\Ber(1+uZ)&={}\sum_{k=0}^{\infty} u^k\ts
\str^{}_{1,\dots,k}\ts A_k Z_1\dots Z_k,\\
\non
\big[\Ber(1-uZ)\big]^{-1}&={}\sum_{k=0}^{\infty}
u^k\ts\str^{}_{1,\dots,k}\ts H_k Z_1\dots Z_k,\\
\non
\big[\Ber(1+uZ)\big]^{-1}\ts\di_u\tss \Ber(1+uZ)
&={}\sum_{k=0}^{\infty} (-u)^k\ts\str\ts Z^{k+1}
\end{align}
relying on the super-analogues of the
MacMahon Master Theorem and Newton identity proved in \cite{mr:mm}.

The center $\z(\wh\g)$ of the vertex algebra
$V_{\text{\rm cri}}(\g)$ is invariant under the {\em translation operator}
$T:V_{\text{\rm cri}}(\g)\to V_{\text{\rm cri}}(\g)$ which is the derivation $T=-d/dt$
of the algebra $\U\big(t^{-1}\g[t^{-1}]\big)$ determined by the properties
\beql{tprop}
\big[T,E_{ij}[r]\big]=-r\tss E_{ij}[r-1],
\eeq
where $E_{ij}[r]$ is understood as the operator of left multiplication by $E_{ij}[r]$.

The following property of the Segal--Sugawara vectors was conjectured
in \cite[Remark~3.4(ii)]{mr:mm}. It
can be regarded as a super-analogue of the Feigin--Frenkel theorem.

\bcj\label{conj:ssgener}
Each of the families $\{T^{\tss r} s_{kk}\}$, $\{T^{\tss r} b_{kk}\}$ and $\{T^{\tss r} h_{kk}\}$
with $r\geqslant 0$ and $k\geqslant 1$
generates the algebra $\z(\wh\g)$.
\qed
\ecj

\subsection{Main results}

We will prove Conjecture~\ref{conj:ssgener} in the case $m=n=1$; that is, for
$\g=\gl(1|1)$. This is
the four-dimensional Lie superalgebra
with the even basis elements $E_{11}, E_{22}$ and odd elements $E_{21}, E_{12}$.
Define the
invariant supersymmetric bilinear form $(.|.)$ on $\g$ by
\ben
(E_{ij}|E_{kl})=\de_{ij}\tss\de_{kl}(-1)^{\bi+\bk}.
\een
By \cite[Sec.~2.5]{k:va},
the corresponding {\em affinization} $\wh\g$ is
the centrally extended Lie superalgebra
of Laurent polynomials
with the commutation relations \eqref{commrel}.
The following is our first main result.

\begin{mthma}
Conjecture~\ref{conj:ssgener} holds for $\g=\gl(1|1)$.
\end{mthma}

As explained in \cite[Chap.~3]{f:lc}
by the example $\g=\sll(2)$, the proof can be reduced to verifying the corresponding property
of the {\em classical limit} of the vacuum
module $V_{\text{\rm cri}}(\g)$; that is, to describing
the invariants of the $\g[t]$-module $\Sr\big(\g[t,t^{-1}]/\g[t]\big)$.
We
regard $\g[t,t^{-1}]/\g[t]$ as a $\g[t]$-module with the adjoint action
and extend it to the symmetric algebra.
We will identify the quotient $\g[t,t^{-1}]/\g[t]$ with
$\wh\g_-=t^{-1}\g[t^{-1}]$ via the natural vector space isomorphism.
The classical limit of $\z(\wh\g)$
is the algebra of invariants
\beql{invs}
\Sr(\wh\g_-)^{\g[t]}
=\{P\in \Sr(\wh\g_-)\ |\ \g[t]\ts P=0\}.
\eeq
In the case where $\g$ is a simple
Lie algebra, the algebra \eqref{invs}
is described by the Beilinson--Drinfeld theorem (see \cite[Theorem~3.4.2]{f:lc}), and
this description is also implied by an earlier work of Ra\"{i}s and Tauvel~\cite{rt:ip}.
Namely, recall that the algebra $\Sr(\g)^{\g}$ of $\g$-invariants in the symmetric algebra $\Sr(\g)$ admits
an algebraically independent family of generators,
\ben
\Sr(\g)^{\g}=\CC[P_1,\dots,P_n],\qquad n=\text{\rm rank }\g;
\een
see e.g. \cite[Sec.~7.3]{d:ae}. Identify the generators $P_i$ with their images
under the embedding $\Sr(\g)\hra \Sr(\wh\g_-)$ taking $X\in\g$ to $X[-1]$.
Then the family $\{T^{\tss r} P_{k}\}$ is algebraically independent and
\ben
\Sr(\wh\g_-)^{\g[t]}=\CC[T^{\tss r} P_1,\dots,T^{\tss r} P_n\ |\ r\geqslant 0],
\een
where $T=-d/dt$ now denotes the derivation
of the algebra $\Sr(\wh\g_-)$ determined by the same properties
\eqref{tprop} applied to the symmetric algebra.

If $\g$ is a simple Lie superalgebra, then the structure of the algebra of invariants
$\Sr(\g)^{\g}$ is more complicated; see e.g. \cite{s:ip}, \cite{s:ir}. In particular,
it does not admit an algebraically independent family of generators.
One could still expect that a natural super-analogue of the
Beilinson--Drinfeld--Ra\"{i}s--Tauvel theorem holds: if $\{P_k\ |\ k\geqslant 1\}$
is a family of generators of $\Sr(\g)^{\g}$, then the derivatives
$\{T^{\tss r}P_k\ |\ k\geqslant 1,\  r\geqslant 0\}$ generate the algebra
$\Sr(\wh\g_-)^{\g[t]}$.
We prove this analogue for $\g=\gl(1|1)$ by producing a basis of the
algebra of invariants and showing that each basis element is expressed
in terms of the generators.
The Chevalley images of the basis elements turn out to form a basis of the algebra
$\La^{\text{\rm aff}}(1|1)$
of {\em affine supersymmetric polynomials}. We establish this property by employing
a lemma of Sergeev~\cite{s:ip, s:ir} which allows us to prove
the second main theorem.

\begin{mthmb}
The algebra of $\g[t]$-invariants
of $\Sr(\wh\g_-)$ is isomorphic to the algebra $\La^{\text{\rm aff}}(1|1)$
of affine supersymmetric polynomials.
\end{mthmb}

In proving Theorem~B, we compare
the Hilbert--Poincar\'{e} series
of these algebras. Both series turn out to coincide
with the generating function of the plane partitions over the $(1,1)$-hook. In more detail,
such a plane partition can be regarded as a finite sequence of Young diagrams
$\la^{(1)}\supset\dots\supset\la^{(r)}$, where each term of the sequence is a hook diagram
$(a,1^b)$. Equivalently, a plane partition can be viewed as a corner ``brick wall" formed
by unit cubes or ``bricks", the $i$-th level of the wall has the shape of the hook $\la^{(i)}$,
as illustrated:


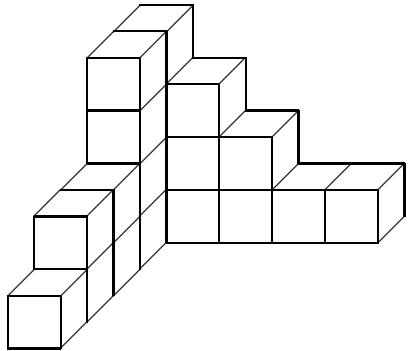
\begin{center}
\begin{picture}(150,130)
\thinlines

\put(0,0){\line(0,1){20}}
\put(10,30){\line(0,1){20}}
\put(20,0){\line(0,1){20}}
\put(30,10){\line(0,1){40}}
\put(30,70){\line(0,1){40}}
\put(40,20){\line(0,1){40}}
\put(50,30){\line(0,1){80}}
\put(60,40){\line(0,1){80}}
\put(70,110){\line(0,1){20}}
\put(80,40){\line(0,1){60}}
\put(90,90){\line(0,1){20}}
\put(100,40){\line(0,1){40}}
\put(110,70){\line(0,1){20}}
\put(120,40){\line(0,1){20}}
\put(140,40){\line(0,1){20}}
\put(150,50){\line(0,1){20}}

\put(0,0){\line(1,0){20}}
\put(0,20){\line(1,0){20}}
\put(10,30){\line(1,0){20}}
\put(10,50){\line(1,0){20}}
\put(20,60){\line(1,0){20}}
\put(30,70){\line(1,0){20}}
\put(30,90){\line(1,0){20}}
\put(30,110){\line(1,0){20}}
\put(40,120){\line(1,0){20}}
\put(50,130){\line(1,0){20}}
\put(60,40){\line(1,0){80}}
\put(60,60){\line(1,0){80}}
\put(60,80){\line(1,0){40}}
\put(60,100){\line(1,0){20}}
\put(70,110){\line(1,0){20}}
\put(90,90){\line(1,0){20}}
\put(110,70){\line(1,0){40}}

\put(0,20){\line(1,1){10}}
\put(10,50){\line(1,1){20}}
\put(20,0){\line(1,1){40}}
\put(20,20){\line(1,1){40}}
\put(30,50){\line(1,1){30}}
\put(30,110){\line(1,1){20}}
\put(50,90){\line(1,1){20}}
\put(50,110){\line(1,1){20}}
\put(80,80){\line(1,1){10}}
\put(80,100){\line(1,1){10}}
\put(100,60){\line(1,1){10}}
\put(100,80){\line(1,1){10}}
\put(120,60){\line(1,1){10}}
\put(140,60){\line(1,1){10}}
\put(140,40){\line(1,1){10}}

\end{picture}
\end{center}

\noindent
The corresponding sequence of hooks in this example is $(5,1^4)\supset (3,1^3)\supset (2,1)\supset (1^2)$.
The generating function of such plane partitions is known \cite{fjmm:qt} and given by the expression
\beql{planep}
\frac{1}{(q)^2_{\infty}}\ts\sum_{k=0}^{\infty}(-1)^k\ts q^{\frac{k^2+k}{2}}
=1+q+3\tss q^2+6\tss q^3+12\tss q^4+21\tss q^5+38\tss q^6+63\tss q^7
+\dots,
\eeq
where
\ben
(q)_{\infty}=\prod_{i=1}^{\infty}(1-q^i)
\een
and the coefficient of $q^N$
is the number of plane partitions with $N$ cubes. We will prove a new
fermionic formula for this generating function.

\begin{mthmc}
The Hilbert--Poincar\'{e} series of the algebra $\La^{\text{\rm aff}}(1|1)$ is given by
\ben
\frac{1}{(q)_{\infty}}\ts\sum_{k=0}^{\infty}\frac{q^{\tss k^2+k}}{(q)^2_k},
\qquad \text{with}\quad (q)_{k}=\prod_{i=1}^{k}(1-q^i).
\een
Moreover, this series coincides with the generating function
of the plane partitions over the $(1,1)$-hook and so equals \eqref{planep}.
\end{mthmc}

We also give a conjectural characterization
property of the affine supersymmetric polynomials
analogous to that of the
supersymmetric polynomials; see Section~\ref{subsec:acp}.
It is implied by the invariance property of elements of the symmetric algebra
in the same way as in the finite-dimensional case; cf. \cite{s:ip}.

In the Appendix we prove a simple formula for the Hilbert--Poincar\'{e}
series of the algebra $\La(m|n)$ of supersymmetric polynomials. In different forms
this series was previously calculated in \cite{oz:sr}
and \cite{sv:dq}.

\medskip

We are grateful to Jonathan Brundan, Vera Serganova and Alexander Sergeev for
useful discussions. The second author is grateful to the University of Sydney
for the hospitality and support during his visit when most of this work
was completed.

\section{Invariants of the symmetric algebra}
\label{sec:invsyma}
\setcounter{equation}{0}

In the specialization $m=n=1$ the matrix $Z=[Z_{ij}]$
takes the form
\ben
Z=\begin{bmatrix}\tau+E_{11}[-1]&E_{12}[-1]\\
-E_{21}[-1]&\tau-E_{22}[-1]
\end{bmatrix}.
\een
Due to the relationship between the elements of the three families
$\{s_{kl}\}$, $\{b_{kl}\}$ and $\{h_{kl}\}$ recalled in the Introduction,
it will be sufficient
to prove Theorem~A for one of them. We will work with the elements $h_{kl}\in\U(\wh\g_-)$
most of the time, whose explicit form implied by \eqref{hkl} is
provided by \cite[Proposition~2.3]{mr:mm}. Denote by $\bar h_{kl}$ their symbols
in the associated graded algebra $\gr\ts\U(\wh\g_-)\cong \Sr(\wh\g_-)$.
They are easily calculated and we have, in particular,
\ben
\bar h_{kk}=E_{11}[-1]^{k-1}\big(E_{11}[-1]+E_{22}[-1]\big)+
(k-1)\tss E_{11}[-1]^{k-2}E_{21}[-1]E_{12}[-1],
\een
where we keep the same notation $E_{ij}[r]$ for the generators of the symmetric
algebra.
Observe that these elements are recovered from the invariants
\ben
E_{11}^{k-1}\big(E_{11}+E_{22}\big)+
(k-1)\tss E_{11}^{k-2}E_{21}E_{12}\in\Sr(\g)^{\g}
\een
through the embedding $\Sr(\g)\hra \Sr(\wh\g_-)$ sending $X\in\g$ to $X[-1]$.
In the same way, the symbols $\bar b_{kk}$ and $\bar s_{kk}$
are obtained as the images in $\Sr(\wh\g_-)$ of the respective $\g$-invariants
\beql{bkkskk}
E_{22}^{k-1}\big(E_{11}+E_{22}\big)-
(k-1)\tss E_{22}^{k-2}E_{21}E_{12}\Fand \str\begin{bmatrix}\phantom{-}E_{11}&\phantom{-}E_{12}\\
-E_{21}&-E_{22}
\end{bmatrix}^k
\eeq
in $\Sr(\g)$.
As before, we regard
the translation operator $T$ as a derivation
on $\Sr(\wh\g_-)$ determined
by \eqref{tprop}.
We have an easily verified generating function identity
\beql{hser}
\sum_{r=0}^{\infty} \frac{T^{\tss r}\tss\bar h_{kk}}{r!}\tss z^r
=E_{11}(z)^{k-1}\tss \big(E_{11}(z)+E_{22}(z)\big)+(k-1)\tss E_{11}(z)^{k-2}\tss E_{21}(z)E_{12}(z),
\eeq
where
\ben
E_{ij}(z)=\sum_{r=0}^{\infty}E_{ij}[-r-1]\tss z^r.
\een
It is a consequence of the commutative
vertex algebra structure on $\Sr(\wh\g_-)$; first we note that
\ben
\sum_{r=0}^{\infty} \frac{T^{\tss r}\tss E_{ij}[-1]}{r!}\tss z^r=E_{ij}(z)
\een
and then apply the property
\ben
\sum_{r=0}^{\infty} \frac{T^{\tss r}\tss \big(E_{ij}[-1]\tss E_{kl}[-1]\big)}{r!}\tss z^r
=E_{ij}(z)\tss E_{kl}(z).
\een
Similarly, the invariants \eqref{bkkskk} give rise to the power series
\beql{bser}
E_{22}(z)^{k-1}\tss \big(E_{11}(z)+E_{22}(z)\big)-(k-1)\tss E_{22}(z)^{k-2}\tss E_{21}(z)E_{12}(z)
\eeq
and
\beql{sser}
\str\begin{bmatrix}\phantom{-}E_{11}(z)&\phantom{-}E_{12}(z)\ts\\
-E_{21}(z)&-E_{22}(z)\ts
\end{bmatrix}^k,
\eeq
respectively.

The next theorem is an analogue of
the Beilinson--Drinfeld--Ra\"{i}s--Tauvel theorem;
see \cite{rt:ip} (and \cite{bb:ei} for a similar argument) and \cite[Theorem~3.4.2]{f:lc}.

\bth\label{thm:bdrt}
The coefficients of each of the series \eqref{hser}, \eqref{bser} and \eqref{sser}
generate the algebra  $\Sr(\wh\g_-)^{\g[t]}$.
\eth

\bpf
The rest of Section~\ref{sec:invsyma} is devoted to the proof of Theorem~\ref{thm:bdrt}.
Our strategy will be to construct a basis of the algebra of
$\g[t]$-invariants of $\Sr(\wh\g_-)$ and then to show that every
basis element can be expressed as a polynomial in the
coefficients of the series \eqref{hser}; i.e., in the
elements $T^{\tss r}\bar h_{kk}$.
It will be convenient to use the following basis elements of
the Lie superalgebra $\wh\g_-$:
for $i\geqslant 0$ set
\beql{ac}
a_i=E_{11}[-i-1],\qquad c_i=E_{11}[-i-1]+E_{22}[-i-1],
\eeq
and
\ben
\vp_i=E_{21}[-i-1],
\qquad \psi_i=E_{12}[-i-1],
\een
so that the $c_i$ are central in $\wh\g_-$.
Moreover, for $i\geqslant 0$
introduce elements of the symmetric algebra $\Sr(\wh\g_-)$ by
\ben
y_i=\sum_{a+b=i}\vp_a\tss\psi_b.
\een

\ble\label{lem:y}
Every element of the algebra of
$\g[t]$-invariants of $\Sr(\wh\g_-)$ can be written as a polynomial
in the elements $a_i, c_i$ and $y_i$ with $i\geqslant 0$.
\ele

\bpf
Any element $P\in \Sr(\wh\g_-)$ can be written in the form
\beql{pexpa}
P=\sum_{I,J}P_{IJ}\tss\vp_{i_1}\dots\vp_{i_n}\psi_{j_1}\dots\psi_{j_m}
\eeq
with the conditions $0\leqslant i_1<\dots<i_n$ and
$j_1>\dots>j_m\geqslant 0$
for uniquely determined polynomials $P_{IJ}$ in the $a_i$ and $c_i$, where
$I=\{i_1,\dots,i_n\}$ and $J=\{j_1,\dots,j_m\}$.
Suppose now that $P$ is invariant.
It will be sufficient to use
only the conditions that $E_{22}[r]\tss P=0$ for $r=0,1$. They are clearly
satisfied by any polynomial in $a_i, c_i$ and $y_i$. For the action of $E_{22}[r]$
we have
\ben
E_{22}[r]\ts \vp_i=\begin{cases} \vp_{i-r}\ \ &\text{if}\ \  i\geqslant r,\\
                                   0\ \ &\text{if}\ \  i< r,
                                   \end{cases}
\qquad\qquad
E_{22}[r]\ts \psi_i=\begin{cases} -\psi_{i-r}\ \ &\text{if}\ \  i\geqslant r,\\
                                   \phantom{-}0\ \ &\text{if}\ \  i< r.
                                   \end{cases}
\een
Hence, $P_{IJ}= 0$ unless $I$ and $J$ have the same cardinality, as implied
by the relation $E_{22}[0]\tss P=0$.
Given such an invariant $P$, we can find an element
\ben
Q=\sum_K Q_K\tss y_{k_1}\dots y_{k_l}\qquad\text{with}\qquad
k_1\geqslant\dots\geqslant k_l\geqslant 0,
\een
where $K=\{k_1,\dots,k_l\}$ and each $Q_K$ is a polynomial in the $a_i$ and $c_i$,
such that the expansion \eqref{pexpa} for $P+Q$
does not contain monomials of the form $\vp_{0}\dots\vp_{n-1}\ts\psi_{j_1}\dots\psi_{j_n}$
for any $n\geqslant 0$.
This follows easily by induction on the $n$-tuples $(j_1,\dots,j_n)$
with the lexicographic order; we assume that if $m<n$ then any $m$-tuple
$(h_1,\dots,h_m)$ precedes
$(j_1,\dots,j_n)$.
Indeed, the largest monomial is eliminated by taking
the sum
\ben
\vp_{0}\dots\vp_{n-1}\ts\psi_{j_1}\dots\psi_{j_n}
+\text{const}\ts y_{j_1}\tss y_{j_2+1} \dots y_{j_n+n-1}
\een
for an appropriate value of the constant.

Furthermore, assuming that none of the monomials
$\vp_{0}\dots\vp_{n-1}\ts\psi_{j_1}\dots\psi_{j_n}$ occurs in $P$,
we will show that $P=0$. Suppose for the contrary that $P\ne 0$ and
take the minimum
$n$-tuple $(i_1,\dots,i_n)$ in the lexicographic order such that
$\vp_{i_1}\dots\vp_{i_n}$ occurs in the expansion of $P$; its coefficient
is a nonzero linear combination of the products $\psi_{j_1}\dots\psi_{j_n}$.
By our assumption, $(i_1,\dots,i_n)=(0,1,\dots,s-1,i_{s+1},\dots,i_n)$
for some $0\leqslant s\leqslant n-1$ and $i_{s+1}>s$.
The condition $E_{22}[1]\tss P=0$ then brings a contradiction since the coefficient
of the monomial $\vp_{0}\vp_1\dots\vp_{s-1}\vp_{i_{s+1}-1}\dots \vp_{i_n}$ in the
expansion of $E_{22}[1]\tss P$ will be nonzero.
\epf

Introduce formal power series
\beql{genfu}
c(z)=\sum_{i=0}^{\infty}c_iz^i,\qquad
\vp(z)=\sum_{i=0}^{\infty}\vp_iz^i,\qquad \psi(z)=\sum_{i=0}^{\infty}\psi_iz^i,\qquad
y(z)=\sum_{i=0}^{\infty}y_iz^i.
\eeq
Note that $y(z)=\vp(z)\tss \psi(z)$ and we have the relations
\beql{relpsiy}
\psi(z)^2=0,\qquad y(z)\tss\psi(z)=0, \qquad y(z)^2=0.
\eeq

Suppose that $z=(z_1,\dots,z_n)$ is a family of independent variables and $\la=(\la_1,\dots,\la_n)$
is a partition.
Its parts are nonnegative integers
satisfying the condition $\la_1\geqslant\dots\geqslant\la_n\geqslant 0$.
The length $\ell(\la)$ is
the number of nonzero parts.
The corresponding Schur polynomial
$s_{\la}(z)$ is defined as the ratio of two alternants,
\ben
s_{\la}(z)=
\frac{\left|\begin{matrix}
z^{\la_1+n-1}_1&\dots&z^{\la_1+n-1}_n\\
\vdots&\vdots&\vdots\\
z^{\la_n}_1&\dots&z^{\la_n}_n
\end{matrix}\right|}{\left|\begin{matrix}
z_1^{n-1}&\dots&z_n^{n-1}\\
\vdots&\vdots&\vdots\\
1&\dots&1
\end{matrix}\right|},
\een
where the denominator is the Vandermonde determinant
\ben
\De=\prod_{i<j}(z_i-z_j);
\een
see e.g. \cite[Ch.~1]{m:sf} for other presentations and properties of the Schur
polynomials.

By the last relation in \eqref{relpsiy},
the power series $y(z_1)\dots y(z_n)$
is divisible by $\De$. The ratio
is skew-symmetric with respect to permutations of the $z_i$ which implies that
$y(z_1)\dots y(z_n)$
is divisible by the square of $\De$.
Since the Schur polynomials form a basis of the algebra of symmetric
polynomials in $z_1,\dots,z_n$,
we can define elements $Y^{(n)}_{\la}\in\Sr(\wh\g_-)$ by the expansion
\beql{expaschur}
\frac{y(z_1)\dots y(z_n)}{\prod_{i\ne j}(z_i-z_j)}=\sum_{\la,\ts \ell(\la)\leqslant n}
Y^{(n)}_{\la} s_{\la}(z).
\eeq

By Lemma~\ref{lem:y}, the algebra of invariants $\Sr(\wh\g_-)^{\g[t]}$ is contained
in the subalgebra $\Sr^{\circ}$ of $\Sr(\wh\g_-)$ generated by the elements
$a_i, c_i$ and $y_i$. The subalgebra of $\Sr^{\circ}$ generated
by the $a_i$ and $c_i$ with $i\geqslant 0$ can be regarded as the algebra of polynomials
in these variables, which we denote by
$\Hr$. We regard $\Sr^{\circ}$ as an $\Hr$-module; elements of $\Hr$ act
by multiplication.

\ble\label{lem:modbas}
The family
\ben
\{Y^{(n)}_{\la}\ |\ n=0,1,\dots,\ \ \ell(\la)\leqslant n\}
\een
forms a basis of the
$\Hr$-module $\Sr^{\circ}$.
\ele

\bpf
The expansion \eqref{expaschur} implies that the family
spans the $\Hr$-module $\Sr^{\circ}$. To prove the linear independence over $\Hr$,
express the elements
$Y^{(n)}_{\la}$ in terms of the generators $\vp_i$ and $\psi_i$.
Since these generators are odd, we have the expansions
\ben
\frac{\vp(z_1)\dots\vp(z_n)}{\De}=\sum_{\mu,\ts \ell(\mu)\leqslant n}\vp_{\mu_1+n-1}\dots\vp_{\mu_n}\tss
s_{\mu}(z)
\een
and
\ben
\frac{\psi(z_1)\dots\psi(z_n)}{\De}=\sum_{\nu,\ts \ell(\nu)\leqslant n}
\psi_{\nu_1+n-1}\dots\psi_{\nu_n}\tss
s_{\nu}(z),
\een
summed over partitions $\mu$ and $\nu$.
Furthermore,
\ben
\frac{y(z_1)\dots y(z_n)}{\prod_{i\ne j}(z_i-z_j)}
=\frac{\vp(z_1)\dots\vp(z_n)\ts\psi(z_1)\dots\psi(z_n)}{\De^2},
\een
and so, taking into account \eqref{expaschur}, we conclude that
\beql{yvppsi}
Y^{(n)}_{\la}=\sum_{\mu,\ts\nu}c^{\la}_{\mu\nu}\ts
\vp_{\mu_1+n-1}\dots\vp_{\mu_n}\ts
\psi_{\nu_1+n-1}\dots\psi_{\nu_n},
\eeq
where the sum is taken over
partitions $\mu$ and $\nu$ of lengths not exceeding $n$,
and the $c^{\la}_{\mu\nu}$ are the Littlewood--Richardson coefficients
defined by the relation
\ben
s_{\mu}(z)s_{\nu}(z)=\sum_{\la} c^{\la}_{\mu\nu}s_{\la}(z),
\een
see e.g. \cite[Ch.~1]{m:sf}. Note that $c^{\la}_{\mu\nu}=0$ unless
$|\la|=|\mu|+|\nu|$, where
$|\la|=\la_1+\dots+\la_n$ denotes the weight of $\la$.
Since $c^{\la}_{\la\varnothing}=1$, and the monomials $\vp_{\la_1+n-1}\dots\vp_{\la_n}\ts
\psi_{n-1}\dots\psi_{0}$ are linearly independent over $\Hr$, then so are the elements
$Y^{(n)}_{\la}$.
\epf

\bre\label{rem:lr}
Two more bases of the $\Hr$-module $\Sr^{\circ}$ (which we will not use) are formed by
the monomials $y_{l_1}\dots y_{l_n}$ with $n\geqslant 0$ and the conditions $l_i-l_{i+1}\geqslant 2$
for $i=1,\dots,n-1$ and $l_n\geqslant 0$ (cf. \cite{sf:fm})
and by the monomials
$y_{k_1}\dots y_{k_n}$ with
$k_1\geqslant\dots\geqslant k_n\geqslant n-1$.
\qed
\ere

Now suppose that $P\in \Sr(\wh\g_-)^{\g[t]}$. By Lemmas~\ref{lem:y} and \ref{lem:modbas},
there is a unique presentation
\beql{pdeco}
P=\sum_{m\geqslant 0}\ \sum_{\mu,\ts\ell(\mu)\leqslant m} P^{(m)}_{\mu}\tss Y^{(m)}_{\mu},
\eeq
where the coefficients $P^{(m)}_{\mu}$ are certain polynomials in $a_i$ and $c_i$.

\ble\label{lem:leadterm}
Given a decomposition \eqref{pdeco} for an invariant $P\in \Sr(\wh\g_-)^{\g[t]}$,
let $n\geqslant 0$ have the property that $P^{(m)}_{\mu}=0$ for all $m>n$
and let a partition $\la=(\la_1,\dots,\la_n)$ be such that
$P^{(n)}_{\la}\ne 0$ but $P^{(n)}_{\mu}=0$ for all $\mu$ with $|\mu|>|\la|$.
Then $P^{(n)}_{\la}$ does not depend on the variables $a_i$ with $i\geqslant n$.
\ele

\bpf
We will use the condition $E_{12}[0]\tss P=0$. The operator $E_{12}[0]$ can be written
in the form
\beql{eonetwo}
E_{12}[0]=-\sum_{j\geqslant 0}\psi_j\di_j+\sum_{r\geqslant 0} c_r\di_{\vp_r},
\eeq
where we denote $\di_j=\di/\di a_j$ and $\di_{\vp_r}$ is the left derivative over $\vp_r$.
The condition on $n$ implies
\ben
\sum_{j\geqslant 0}\psi_j\sum_{\mu,\ts\ell(\mu)\leqslant n}
\di_j(P^{(n)}_{\mu})\ts Y^{(n)}_{\mu}=0.
\een
Take $i\geqslant n$ and consider the coefficient of the monomial
$\vp_{\la_1+n-1}\dots\vp_{\la_n}\ts
\psi_i\tss\psi_{n-1}\dots\psi_0$ on the left hand side.
By the condition on $\la$, this monomial can only occur for $j=i$ and $\mu=\la$
thus implying $\di_i(P^{(n)}_{\la})=0$, as required.
\epf

In what follows we will call by a {\em leading component}
any product of the form $R\ts Y^{(n)}_{\la}$,
where $R$ is a polynomial in the $a_i$ and $c_i$ which does
not depend on the variables $a_i$ with $i\geqslant n$.
By Lemma~\ref{lem:leadterm}, every invariant has a leading component.
Our next goal is to show that there exists an invariant $P\in \Sr(\wh\g_-)^{\g[t]}$
containing any given leading component and no other
leading components in the
expansion \eqref{pdeco}. It suffices to do this for monomials
of the form
\beql{ynmon}
Y^{(n)}_{\la}\ts\di_0^{-k_0}\dots\di_{n-1}^{-k_{n-1}}\ts 1,
\eeq
where we regard $\di_i^{-1}$ as a partial integration operator with respect to $a_i$
so that
\beql{amon}
\di_0^{-k_0}\dots\di_{n-1}^{-k_{n-1}}\ts 1=
\frac{a_0^{k_0}\dots a_{n-1}^{k_{n-1}}}{k_0!\dots k_{n-1}!}.
\eeq

Given a value of $n$, introduce another family of independent
variables $t_0,\dots,t_{n-1}$ and
consider the formal power series in the variables $z_i$ and $t_i$
whose coefficients are polynomials in the $a_i$,
\begin{multline}
F(z_1,\dots,z_n;t_0,\dots,t_{n-1})=\prod_{i=0}^{n-1}\Big(1-\di_i^{-1}t_i\Big)^{-1}\\
{}\times\prod_{j=0}^{\infty}\Bigg(1-\di_{n+j}^{-1}\Big(t^{}_{n-1}\tss s^{}_{(j+1)}(z)
-t^{}_{n-2}\tss s^{}_{(j+1,1)}(z)
+\dots+(-1)^{n-1}t^{}_0\tss s^{}_{(j+1,1^{n-1})}(z)\Big)\Bigg)^{-1}\ts 1,
\non
\end{multline}
where $s_{(j+1,1^{k})}(z)$ is the Schur polynomial in the variables $z_1,\dots,z_n$
associated with the hook partition $(j+1,1^{k})$. In particular,
the series $F(z_1,\dots,z_n;t_0,\dots,t_{n-1})$ is symmetric in $z_1,\dots,z_n$.

\ble\label{lem:formf}
We have the identity
\ben
F(z_1,\dots,z_n;t_0,\dots,t_{n-1})=\prod_{i=0}^{\infty}\Bigg(1-\di_i^{-1}
\Big(z_1^i\tss T_n^{(1)}+\dots+z_n^i\tss T_n^{(n)}\Big)\Bigg)^{-1}\ts 1,
\een
where
\ben
T_n^{(k)}=\frac{t_{n-1}-t_{n-2}\ts e_1(z_1,\dots,\wh z_k,\dots,z_n)+\dots
+(-1)^{n-1}t_0\ts e_{n-1}(z_1,\dots,\wh z_k,\dots,z_n)}
{(z_k-z_1)\ldots\wedge\ldots(z_k-z_n)}
\een
and $e_1,\dots,e_{n-1}$ denote the elementary symmetric polynomials;
the hats and wedges indicate symbols or zero factors to be skipped.
\ele

\bpf
The rational function $T_n^{(k)}$ is written as the ratio
\ben
T_n^{(k)}
=\frac{\left|\begin{matrix}
z_1^{n-1}&\dots&t_{n-1}&\dots&z_n^{n-1}\\
\vdots&\vdots&\vdots&\vdots&\vdots\\
z_1&\dots&t_1&\dots&z_n\\
1&\dots&t_0&\dots&1
\end{matrix}\right|}{\De},
\een
where the $t_i$ occupy the $k$-th column in the numerator. Hence, the $T_n^{(k)}$
are the solutions of the system of equations
\ben
z_1^i\tss T_n^{(1)}+\dots+z_n^i\tss T_n^{(n)}=t_i,\qquad i=0,1,\dots,n-1.
\een
Furthermore, if $i\geqslant n$ and $1\leqslant m\leqslant n$
then the coefficient of $t_{n-m}$
in the expression $z_1^i\tss T_n^{(1)}+\dots+z_n^i\tss T_n^{(n)}$ equals
the ratio
\ben
\frac{\left|\begin{matrix}
z_1^{n-1}&\dots&z_n^{n-1}\\
\dots&\dots&\dots\\
z_1^i&\dots&z_n^i\\
\dots&\dots&\dots\\
1&\dots&1
\end{matrix}\right|}{\De},
\een
where $z_1^i,\dots,z_n^i$ replace row $m$ of the Vandermonde determinant
in the numerator. This ratio coincides with $(-1)^{m-1}\tss s_{(i-n+1,1^{m-1})}$,
as required.
\epf

We are now in a position to prove a key lemma providing explicit $\g[t]$-invariants
in $\Sr(\wh\g_-)$. Recall the formal power series \eqref{genfu} and set
\beql{az}
A(z_1,\dots,z_n;t_0,\dots,t_{n-1})=\prod_{k=1}^n
\Big(c(z_k)+y(z_k)\tss T_n^{(k)}\Big)
F(z_1,\dots,z_n;t_0,\dots,t_{n-1}).
\eeq
This is a formal power series in the $z_i$ and $t_i$, symmetric in $z_1,\dots,z_n$,
whose coefficients are
elements of the subalgebra $\Sr^{\circ}$ of $\Sr(\wh\g_-)$ generated by the $a_i, c_i$ and $y_i$.

\ble\label{lem:inv}
All coefficients of the series $A(z_1,\dots,z_n;t_0,\dots,t_{n-1})$ belong
to $\Sr(\wh\g_-)^{\g[t]}$.
\ele

\bpf
It is enough to
show that $E_{12}[0]\tss A(z_1,\dots,z_n;t_0,\dots,t_{n-1})=0$.
Recall that the action of $E_{12}[0]$ is given by the operator
\eqref{eonetwo}.
Lemma~\ref{lem:formf} implies
\ben
\di_i\tss F(z_1,\dots,z_n;t_0,\dots,t_{n-1})=\Big(z_1^i\tss T_n^{(1)}+
\dots+z_n^i\tss T_n^{(n)}\Big)F(z_1,\dots,z_n;t_0,\dots,t_{n-1})
\een
for all $i\geqslant 0$.
Hence,
\ben
E_{12}[0]\tss F(z_1,\dots,z_n;t_0,\dots,t_{n-1})=
-\Big(\psi(z_1)\tss T_n^{(1)}+
\dots+\psi(z_n)\tss T_n^{(n)}\Big)
\tss F(z_1,\dots,z_n;t_0,\dots,t_{n-1}).
\een
On the other hand, $E_{12}[0]\tss y(z_k)=c(z_k)\tss \psi(z_k)$ and so
\ben
E_{12}[0]\tss\prod_{k=1}^n
\Big(c(z_k)+y(z_k)\tss T_n^{(k)}\Big)=\sum_{i=1}^n
c(z_i)\tss \psi(z_i)\tss T_n^{(i)}\ts\prod_{k\ne i}
\Big(c(z_k)+y(z_k)\tss T_n^{(k)}\Big).
\een
Since $y(z_k)\psi(z_k)=0$ by \eqref{relpsiy}, we have
$E_{12}[0]\tss A(z_1,\dots,z_n;t_0,\dots,t_{n-1})=0$.
\epf

Now expand \eqref{az} along the basis formed by the products of monomials in the $t_i$
and Schur polynomials in $z_1,\dots,z_n$.
Take a partition $\la$ of length not exceeding $n$
and consider the coefficient of the basis element
\beql{basts}
t_0^{k_0}\dots t_{n-2}^{k_{n-2}}
t_{n-1}^{k_{n-1}+n}\tss s_{\la}(z),\qquad k_i\geqslant 0,
\eeq
in the expansion. Furthermore, use Lemma~\ref{lem:modbas} to write this
coefficient as a linear combination of the basis elements $Y^{(m)}_{\mu}$.
By \eqref{expaschur}, this linear combination contains
a leading component in the form \eqref{ynmon}. All other elements $Y^{(m)}_{\mu}$
occurring in the linear combination will have the property $m\leqslant n$; moreover,
if $m=n$ then $|\mu|<|\la|$. Therefore, eliminating all other leading components
with the use of an easy induction, we get an invariant containing a unique leading
component.
Thus, taking into account Lemma~\ref{lem:leadterm},
we may conclude that the coefficients of the basis elements \eqref{basts}
in the expansion of \eqref{az}
with $n$ running over nonnegative integers form a basis of $\Sr(\wh\g_-)^{\g[t]}$
as a module over the algebra of polynomials in the $c_i$ with $i\geqslant 0$.

Take $n=1$ in \eqref{az} and observe that the coefficient of $t_0^{k-1}$
in $A(z_1;t_0)$ equals
\beql{aztone}
\frac{1}{(k-1)!}\Big(a(z)^{k-1}\tss c(z)+(k-1)\tss a(z)^{k-2}\tss y(z)\Big),
\eeq
where
\ben
a(z)=\sum_{i=0}^{\infty} a_iz^i.
\een
This is immediate from the identity
\ben
\sum_{0\leqslant i_1\leqslant\dots\leqslant i_p}z^{i_1+\dots+i_p}\tss
\di^{-1}_{i_1}\dots \di^{-1}_{i_p}\ts 1=\frac{a(z)^p}{p!}
\een
which holds for any $p\geqslant 0$. Since the series \eqref{hser} equals $(k-1)!$ times
\eqref{aztone}, the proof of Theorem~\ref{thm:bdrt}
will be completed if we show that all coefficients of the series \eqref{az}
for all values of $n$ are expressed as polynomials in the coefficients of $A(z_1;t_0)$.
This is the statement of the next lemma.

\ble\label{lem:bd}
We have the identity
\ben
A(z_1,\dots,z_n;t_0,\dots,t_{n-1})=
A(z_1;T_n^{(1)})\dots A(z_n;T_n^{(n)}).
\een
\ele

\bpf
We have
\beql{azt}
A(z_1;t_0)\dots A(z_n;t_{n-1})=\prod_{k=1}^n
\Big(c(z_k)+y(z_k)\tss t_{k-1}\Big)\ts F(z_1;t_0)\dots F(z_n;t_{n-1}).
\eeq
Write
\ben
F(z_1;t_0)\dots F(z_n;t_{n-1})=\prod_{i=0}^{\infty}\Big(1-\di_i^{-1}
z_1^i\tss t_0\Big)^{-1}\tss 1\dots \Big(1-\di_i^{-1}
z_n^i\tss t_{n-1}\Big)^{-1}\tss 1.
\een
Expanding the series and using the identity
\ben
\di_i^{-k_1}\tss 1\dots \di_i^{-k_n}\tss 1=\binom{k_1+\dots+k_n}{k_1,\dots,k_n}\tss
\di_i^{-k_1-\dots-k_n}\tss 1
\een
we find that
\ben
F(z_1;t_0)\dots F(z_n;t_{n-1})=\prod_{i=0}^{\infty}\Bigg(1-\di_i^{-1}
\Big(z_1^i\tss t_0+\dots+z_n^i\tss t_{n-1}\Big)\Bigg)^{-1}\tss 1.
\een
Hence, replacing $t_i\mapsto T_n^{(i+1)}$ for $i=0,\dots,n-1$ in \eqref{azt}
we recover the formal power series $A(z_1,\dots,z_n;t_0,\dots,t_{n-1})$,
as required.
\epf
This completes the proof of Theorem~\ref{thm:bdrt}.
\epf

\section{Affine supersymmetric polynomials}
\label{sec:asp}
\setcounter{equation}{0}

Recall that a polynomial $P(u,v)=P(u_1,\dots,u_m,v_1,\dots,v_n)$
in two sets of independent variables $u=(u_1,\dots,u_m)$
and $v=(v_1,\dots,v_n)$ is called {\em supersymmetric}, if it is symmetric
in each of the sets separately and the following cancellation property holds:
the result of the substitution $u_m=-v_n=t$ into $P(u,v)$ is independent of $t$.
We denote the algebra of supersymmetric polynomials by $\La(m|n)$. The supersymmetric
Schur polynomials parameterized by all Young diagrams not containing the box $(m+1,n+1)$
form a basis of this algebra. Moreover, each of the families of elementary,
complete and power sums supersymmetric
functions generates $\La(m|n)$; see e.g. \cite[Ch.~1]{m:sf}.

Given a supersymmetric polynomial $P(u,v)$, replace each variable $u_i$ and $v_j$ by
the respective formal power series
\ben
u_i(z)=\sum_{r=0}^{\infty}u_{i\tss r}z^r,\qquad v_j(z)=\sum_{r=0}^{\infty}v_{j\tss r}z^r,
\een
and write
\ben
P\big(u_1(z),\dots,u_m(z),v_1(z),\dots,v_n(z)\big)=\sum_{r=0}^{\infty} P_r\tss z^r,
\een
where the coefficients $P_r$ are polynomials in the variables
$u_{1\tss r},\dots,u_{m\tss r},v_{1\tss r},\dots,v_{n\tss r}$
with $r$ running over the set of nonnegative integers.
Equivalently, $P_r$ is found as the derivative
\ben
P_r=\frac{T^r P}{r!},
\een
where $P=P(u,v)$ is regarded as a polynomial in the variables $u_{i\tss 0}=u_i$
and $v_{j\tss 0}=v_j$, and the derivation $T$ acts on the variables by
the rule (cf. Sec.~\ref{sec:invsyma}):
\ben
T:u_{i\tss r}\mapsto (r+1)\tss u_{i\tss r+1},\qquad
v_{j\tss r}\mapsto (r+1)\tss v_{j\tss r+1}.
\een

\bde\label{def:asfser}
We denote by $\La^{\text{\rm aff}}(m|n)$ the subalgebra
of the algebra of polynomials in the variables $u_{i\tss r}$ and $v_{j\tss r}$
generated by all coefficients $P_r$ associated with all supersymmetric
polynomials $P(u,v)$. Any element of $\La^{\text{\rm aff}}(m|n)$ will be called
an {\em affine supersymmetric polynomial}.
\qed
\ede

It is clear that the algebra $\La^{\text{\rm aff}}(m|n)$ is generated by the
coefficients $P_r$ associated to any family $\{P\}$ of generators of the algebra $\La(m|n)$.
For instance, considering the supersymmetric power sums
\ben
u_1^k+\dots+u_m^k-(-1)^{k}(v_1^k+\dots+v_n^k)
\een
we get the following explicit
formulas for generators of $\La^{\text{\rm aff}}(m|n)$:
\beql{genkr}
\sum_{i=1}^m\sum_{r_1+\dots+r_k=r}u_{i\tss r_1}\dots u_{i\tss r_k}
-(-1)^{k}\ts\sum_{j=1}^n\sum_{r_1+\dots+r_k=r}v_{j\tss r_1}\dots v_{j\tss r_k},\qquad
k\geqslant 1,\quad r\geqslant 0,
\eeq
where the second sums are taken over the $k$-tuples $(r_1,\dots,r_k)$
of nonnegative integers.

Setting
\ben
\deg u_{i\tss r}=r+1\Fand  \deg v_{j\tss r}=r+1
\een
defines a grading on the algebra of polynomials in the $u_{i\tss r}$ and $v_{j\tss r}$.
In particular, the degree of the generator in \eqref{genkr} equals $k+r$.
The subalgebra $\La^{\text{\rm aff}}(m|n)$ inherits the grading so that we have the direct
sum decomposition
\ben
\La^{\text{\rm aff}}(m|n)=\bigoplus_{N\geqslant 0} \La^{\text{\rm aff}}(m|n)^N,
\een
where $\La^{\text{\rm aff}}(m|n)^N$ denotes the subspace of $\La^{\text{\rm aff}}(m|n)$ spanned
by homogeneous elements of degree $N$ and we set $\La^{\text{\rm aff}}(m|n)^0:=\CC$.
We let $H_{m,n}(q)$ denote the corresponding Hilbert--Poincar\'{e} series
\ben
H_{m,n}(q)=\sum_{N=0}^{\infty} \dim \La^{\text{\rm aff}}(m|n)^N\ts q^N.
\een

As in the Introduction, by a plane partition over the $(m,n)$-hook we mean
a finite sequence of Young diagrams (or partitions)
$\la^{(1)}\supset\dots\supset\la^{(r)}$ such that $\la^{(1)}$ does not contain
the box $(m+1,n+1)$. Such a plane partition can be viewed as an array formed
by unit cubes, the $i$-th level of the array has the shape $\la^{(i)}$.
An explicit formula for the generating function of the plane
partitions was conjectured in \cite{fjmm:qt}
and proved in \cite{fm:qt}.
For $n\geqslant m\geqslant 1$ it has the form
\ben
\bal
f_{m,n}(q)=\frac{1}{(q)^{m+n}_{\infty}}\ts\sum_{k_1\geqslant\dots\geqslant k_m\geqslant 0}
&\Bigg((-1)^{k_1+\dots+k_m}\ts q^{\frac12\tss \sum_{i=1}^m(k_i^2+(2i-1)k_i)}\\
{}&\times{} \prod_{1\leqslant i<j\leqslant m}(1-q^{k_i-k_j+j-i})
\prod_{1\leqslant i<j\leqslant n}(1-q^{k_i-k_j+j-i})\Bigg),
\eal
\een
where $k_j:=0$ for $j>m$ and the coefficient of $q^N$ in the series is the number
of plane partitions over the $(m,n)$-hook containing exactly $N$ unit cubes.

\bcj\label{conj:hpseries}
The dimension $\dim \La^{\text{\rm aff}}(m|n)^N$ equals the number of plane
partitions over the $(m,n)$-hook containing exactly $N$ unit cubes.
Equivalently, if $n\geqslant m\geqslant 1$ then the Hilbert--Poincar\'{e} series
$H_{m,n}(q)$ coincides with $f_{m,n}(q)$.
\qed
\ecj

The conjecture holds for $n=0$ (or $m=0$); that is, for the algebra of {\em affine symmetric
polynomials} $\La^{\text{\rm aff}}(m)$. This algebra admits a family of algebraically
independent generators which can be obtained, for instance, by taking $n=0$ in \eqref{genkr}:
\ben
\sum_{i=1}^m\sum_{r_1+\dots+r_k=r}u_{i\tss r_1}\dots u_{i\tss r_k},\qquad k=1,\dots,m,
\qquad r\geqslant 0.
\een
The Hilbert--Poincar\'{e} series is then found by
\ben
\prod_{k=1}^m\ts\prod_{r\geqslant k}\ts(1-q^{r})^{-1}=\frac{1}{(q)^{m}_{\infty}}\ts
\prod_{i=1}^{m-1}\ts(1-q^{i})^{m-i}
\een
which coincides with $f_{0,m}(q)$; cf. \cite[Sec.~4.3]{f:lc}.

Below we prove Conjecture~\ref{conj:hpseries} for $m=n=1$; see Sec.~\ref{sec:ch}.
First we give an alternative expression
for the generating function $f_{m,n}(q)$ in this case.

\bpr\label{prop:fermi}
We have
\ben
f_{1,1}(q)=\frac{1}{(q)_{\infty}}\ts\sum_{k=0}^{\infty}\frac{q^{\tss k^2+k}}{(q)^2_k}.
\een
\epr

\bpf
By definition,
\ben
f_{1,1}(q)=\frac{1}{(q)^2_{\infty}}\ts\sum_{k=0}^{\infty}(-1)^k\ts q^{\frac{k^2+k}{2}}.
\een
The desired identity follows from a more general relation which holds for
$s\geqslant 0$:
\beql{ids}
\sum_{k=0}^{\infty}\frac{q^{\tss k^2+k}}{(q)^2_k}-
\frac{1}{(q)_{\infty}}\ts\sum_{k=0}^{s-1}(-1)^k\ts q^{\frac{k^2+k}{2}}
=(-1)^s\ts \sum_{k=s}^{\infty}\frac{q^{\tss k^2-(s-1)k+\frac{s^2-s}{2}}}{(q)_k\ts (q)_{k-s}}.
\eeq
We prove \eqref{ids} by induction on $s$. It holds trivially for $s=0$ so suppose that
$s\geqslant 1$. To complete the induction step we need to show that
\ben
(-1)^s\ts \sum_{k=s}^{\infty}\frac{q^{\tss k^2-(s-1)k+\frac{s^2-s}{2}}}{(q)_k\ts (q)_{k-s}}
-(-1)^s\ts \frac{q^{\frac{s^2+s}{2}}}{(q)_{\infty}}
=(-1)^{s+1}\ts \sum_{k=s+1}^{\infty}\frac{q^{\tss k^2-s\tss k+\frac{s^2+s}{2}}}{(q)_k\ts (q)_{k-s-1}}.
\een
This is immediate from the identity
\ben
\frac{1}{(q)_{\infty}}=\sum_{k=s}^{\infty}\frac{q^{\tss k(k-s)}}{(q)_k\ts (q)_{k-s}}
\een
which holds for $s\geqslant 0$ and is easily verified as follows.
Both sides are generating functions for
all partitions. This is clear for the left hand side, while the expression
on the right hand side is obtained by first assigning the maximum size rectangle of the form
$(k-s)\times k$ contained in a Young diagram. Then the generating function of the Young
diagrams with a fixed value of $s$ is given by
\ben
\frac{q^{\tss k(k-s)}}{(q)_k\ts (q)_{k-s}}
\een
as required.
\epf

\subsection{Affine cancellation property}
\label{subsec:acp}

Using notation \eqref{ac}, we will regard $\La^{\text{\rm aff}}(1|1)$
as the subalgebra of the algebra of polynomials in the variables $a_r=u_{1\tss r}$
and $c_r=u_{1\tss r}+v_{1\tss r}$ with $r\geqslant 0$.
Working over Laurent polynomials in $c_0$ define elements $d_r$ by the relation
\ben
d(z):=\sum_{r=0}^{\infty} d_rz^r=c(z)^{-1},\qquad c(z)=\sum_{r=0}^{\infty} c_rz^r.
\een
Explicitly,
\ben
d_r=c_0^{-1}\ts\sum_{\al_1+2\ts\al_2+\dots+r\ts\al_r=r}
\frac{(\al_1+\dots+\al_r)!}{\al_1!\ts\al_2!\dots\al_r!}\ts
\Big({-}\frac{c_1}{c_0}\Big)^{\al_1}\dots \Big({-}\frac{c_r}{c_0}\Big)^{\al_r},
\een
summed over nonnegative integers $\al_i$.
Consider the operator
\ben
\Dc=\sum_{r=0}^{\infty} d_r\ts\di_r,\qquad \di_r=\di/\di\tss a_r.
\een

\bpr\label{prop:canc}
If $P\in \La^{\text{\rm aff}}(1|1)$ then $\Dc\tss P$
does not contain negative powers of $c_0$.
\epr

\bpf
The algebra $\La^{\text{\rm aff}}(1|1)$ is generated by the coefficients
of the series $a(z)^{k}\tss c(z)$ with $k\geqslant 0$.
We have
\ben
\Dc\tss a(z)^{k}\tss c(z)=k\tss a(z)^{k-1}\tss d(z)\tss c(z)=k\tss a(z)^{k-1}.
\een
Thus, the required property holds for generators of the algebra $\La^{\text{\rm aff}}(1|1)$.
Since $\Dc$ is a derivation, it will hold for
all its elements.
\epf

We conjecture that the property given by Proposition~\ref{prop:canc}
is characteristic for the affine supersymmetric polynomials.

\bcj\label{conj:canc}
A polynomial $P$ in the variables $a_r$ and $c_r$ belongs to $\La^{\text{\rm aff}}(1|1)$
if and only if $\Dc\tss P$
does not contain negative powers of $c_0$.
\ecj

\section{Chevalley-type isomorphism}
\label{sec:ch}
\setcounter{equation}{0}

In this section we prove Theorems B and C.

Let $\g=\n_-\oplus\h\oplus\n_+$ be the triangular decomposition of $\g=\gl(m|n)$,
where the subalgebras $\n_-,\h$ and $\n_+$ are spanned
by the basis elements $E_{ij}$ with $i<j$, $i=j$ and $i>j$, respectively.
The Chevalley homomorphism
\ben
\vs:\Sr(\g)\to \Sr(\h)
\een
is the projection modulo the ideal $\Sr(\g)(\n_-\cup\n_+)$. The restriction
of $\vs$ to the subalgebra of invariants yields an isomorphism between $\Sr(\g)^{\g}$
and the algebra of supersymmetric polynomials in two sets
of variables; see
e.g. \cite{s:ip}.

Consider an affine analogue of $\vs$ defined as the projection
\beql{affche}
\wh\vs:\Sr(\wh\g_-)\to \Sr(\wh\h_-)
\eeq
modulo the ideal $\Sr(\wh\g_-)\big(t^{-1}\n_-[t^{-1}]\cup t^{-1}\n_+[t^{-1}]\big)$,
where we set $\wh\h_-=t^{-1}\h[t^{-1}]$. We identify $\Sr(\wh\h_-)$ with the algebra
of polynomials in the variables $u_{1\tss r},\dots,u_{m\tss r},v_{1\tss r},\dots,v_{n\tss r}$
with $r\geqslant 0$ by setting
\ben
u_{i\tss r}=E_{i\ts i}[-r-1]\Fand v_{j\tss r}=E_{j+m\ts j+m}[-r-1].
\een

\bpr\label{prop:chev}
The restriction of the homomorphism \eqref{affche} to the subalgebra
$\Sr(\wh\g_-)^{\g[t]}$ is injective.
\epr

\bpf
Suppose that $Q\in \Sr(\wh\g_-)^{\g[t]}$ and $\wh\vs\tss(Q)=0$.
Take a positive integer $p$ such that $Q$ does not depend on the generators
$E_{ij}[r]$ with $r< -p$. By the definition of the $\g[t]$-action on
$\Sr(\wh\g_-)$, we have $t^p\tss\g[t]\ts Q=0$. Denote by $\g_p$ the quotient of $\g[t]$
by the ideal $t^p\tss\g[t]$ and denote by $\g_{p,-}$
the quotient of $\wh\g_-=t^{-1}\g[t^{-1}]$ by the ideal $t^{-p-1}\g[t^{-1}]$.
The proposition will follow if we show that the restriction map
\ben
\Sr(\g_{p,-})^{\g_p}\to \Sr(\h_{p,-})
\een
is injective for any positive integer $p$, where $\h_{p,-}$ denotes the quotient
of $\wh\h_-$ by the ideal $t^{-p-1}\h[t^{-1}]$.
We will derive this claim from the following general
result of Sergeev~\cite[Proposition~1.1]{s:ip};
see also \cite[Lemma~4.3]{s:ir} for a shorter and more direct proof.

\ble\label{lem:ser}
Let $\g$ be a finite-dimensional Lie superalgebra and $V$
a finite-dimensional $\g$-module. Given a subspace $W$ of $V$, suppose that
there exists an even element $w_0\in W$ such that the map
\ben
\g\times W\to V,\qquad (x,w)\mapsto x\tss w_0+w
\een
is surjective. Then the restriction map $\Sr(V^*)^{\g}\to\Sr(W^*)$
is injective.
\qed
\ele

To apply the lemma we take
$\g=\g_p$ and let $V=\g_p$ be the adjoint $\g_p$-module. The dual module $V^*$
is isomorphic to $\g_{p,-}$; the isomorphism takes the element $E_{kl}[s]^*$ dual to
the basis vector $E_{kl}[s]$ of $\g_p$ to the element $E_{lk}[-s-1](-1)^{\bar l}$.
The subspace $W$ is the quotient $\h_p$ of $\h[t]$ by the ideal $t^p\tss \h[t]$.
The dual space $W^*$ is identified with $\h_{p,-}$. The assumptions of Lemma~\ref{lem:ser}
will hold for any element
\ben
w_0=\sum_{k=1}^{m+n}\ga_k\tss E_{kk}[0]
\een
with $\ga_i\ne \ga_j$ for $i\ne j$. Indeed, this is clear from the relations
\ben
[E_{ij}[r],w_0]=(\ga_j-\ga_i)\ts E_{ij}[r].
\een
The proposition is proved.
\epf

Now we return to the Lie superalgebra
$\g=\gl(1|1)$. The grading on the symmetric algebra $\Sr(\wh\g_-)$
is defined by
\ben
\deg E_{ij}[-r-1]=r+1.
\een
We have $\deg \vp_i=\deg \psi_i=i+1$, and by \eqref{yvppsi}
\beql{degy}
\deg Y^{(k)}_{\la}=|\la|+k\tss(k+1),\qquad \ell(\la)\leqslant k.
\eeq
For a given $k$, let $d_N$ be the number of basis elements $Y^{(k)}_{\la}$
of degree $N$. By \eqref{degy},
the generating function is given by
\ben
\sum_{N\geqslant 0}d_N\ts q^N=\frac{q^{\tss k^2+k}}{(q)_k}.
\een
Since $\deg a_i=i+1$ the generating function of the monomials \eqref{amon}
(with $n:=k$)
is $(q)^{-1}_k$. Similarly, the generating function of the algebra
of polynomials in the $c_i$ is $(q)^{-1}_{\infty}$ and so
the Hilbert--Poincar\'{e} series of $\Sr(\wh\g_-)^{\g[t]}$ is given by
\ben
\frac{1}{(q)_{\infty}}\ts\sum_{k=0}^{\infty}\frac{q^{\tss k^2+k}}{(q)^2_k}.
\een
Furthermore, the image of the series \eqref{hser}
under the projection \eqref{affche} equals
\beql{seon}
E_{11}(z)^{k-1}\tss \big(E_{11}(z)+E_{22}(z)\big).
\eeq
Since the elements $u^{k-1}_1\tss(u_1+v_1)$ with $k\geqslant 1$ generate the
algebra of supersymmetric polynomials $\La(1|1)$, the coefficients
of the series \eqref{seon} generate the algebra
$\La^{\text{\rm aff}}(1|1)$ of affine supersymmetric polynomials. Here, as before,
we identify the variables by $u_{1\tss r}=E_{11}[-r-1]$ and $v_{1\tss r}=E_{22}[-r-1]$.
Therefore, by Proposition~\ref{prop:chev}, we have a Chevalley-type isomorphism
of graded algebras
\ben
\Sr(\wh\g_-)^{\g[t]}\cong \La^{\text{\rm aff}}(1|1).
\een
This completes the proof of Theorem B while Theorem C now follows
from Proposition~\ref{prop:fermi}.

\bcj\label{conj:chevalley}
Let $\g=\gl(m|n)$.
The restriction of the map \eqref{affche} to the subalgebra of $\g[t]$-invariants
yields an isomorphism of graded algebras
\ben
\Sr(\wh\g_-)^{\g[t]}\cong \La^{\text{\rm aff}}(m|n).
\een
In particular, their Hilbert--Poincar\'{e} series coincide.
\qed
\ecj

Besides the case $m=n=1$,
Conjecture~\ref{conj:chevalley} also holds for $n=0$ or $m=0$ as implied by the
Beilinson--Drinfeld--Ra\"{i}s--Tauvel theorem; see \cite[Sec.~4.3]{f:lc}.

\appendix

\section{Generating function for the supersymmetric polynomials
in $m+n$ variables}
\label{sec:genf}
\setcounter{section}{1}
\setcounter{equation}{0}

Two different forms of the Hilbert--Poincar\'{e} series
for the algebra $\La(m|n)$ were given in \cite{oz:sr} and \cite{sv:dq}.
Both proofs rely on the parametrization of basis elements
of $\La(m|n)$ by Young diagrams contained in the $(m,n)$-hook.
We give yet another formula for the series and derive it from
the characterization of the supersymmetric polynomials via the cancellation property;
cf. \cite{s:cs}.

\bpr\label{prop:hpser}
The Hilbert--Poincar\'{e} series of the algebra $\La(m|n)$ is found by
\ben
\chi_{m,n}(q)=\sum_{k=0}^{\min\{m,n\}}\frac{q^{(m-k)(n-k)}}{(q)_{m-k}\ts (q)_{n-k}}.
\een
\epr

\bpf
As before, we consider two sets of variables $u=(u_1,\dots,u_m)$
and $v=(v_1,\dots,v_n)$. For $m,n\geqslant 1$
we have a surjective homomorphism
\ben
\La(m,n)\to \La(m-1,n-1),
\qquad u_m\mapsto 0,\quad v_n\mapsto 0.
\een
Its kernel coincides with the space
\ben
\prod_{i=1}^m\tss\prod_{j=1}^n \ts(u_i+v_j)\ts \CC[u,v]^{\Sym_m\times\Sym_n},
\een
where $\CC[u,v]^{\Sym_m\times\Sym_n}$ is the algebra of bisymmetric
polynomials. Hence we have a recurrence relation
\ben
\chi_{m,n}(q)=\chi_{m-1,n-1}(q)+\frac{q^{mn}}{(q)_{m}\ts (q)_{n}}
\een
which leads to the desired formula.
\epf

The recurrence relation can also be easily seen from
the parametrization of basis elements
by Young diagrams. The term $\chi_{m-1,n-1}(q)$ accounts for the diagrams
contained in the $(m-1,n-1)$-hook, whereas the generating function of
the diagrams in the $(m,n)$-hook containing the box $(m,n)$ is
$q^{mn}/(q)_{m}\ts (q)_{n}$; cf. the proof of Proposition~\ref{prop:fermi}.


\begin{thebibliography}{99}


\bibitem{bb:ei}
J. Brown and J. Brundan,
{\it Elementary invariants for centralizers
of nilpotent matrices}, J. Aust. Math. Soc. {\bf 86} (2009), 1--15.

\bibitem{cm:ho}
A. V. Chervov and A. I. Molev,
{\it On higher order Sugawara operators},
Int. Math. Res. Not. (2009), no. 9, 1612--1635.

\bibitem{ct:qs}
A. Chervov and D. Talalaev,
{\it Quantum spectral curves,
quantum integrable systems  and
the geometric Langlands correspondence},
{\tt arXiv:hep-th/0604128}.

\bibitem{d:ae}
{J. Dixmier},
{\it Alg\`ebres Enveloppantes},
{Gauthier-Villars, Paris},
1974.

\bibitem{ff:ak}
B. Feigin and E. Frenkel,
{\it Affine Kac--Moody algebras at the critical level
and Gelfand--Dikii algebras},
Int. J. Mod. Phys. A{\bf 7}, Suppl. 1A (1992), 197--215.

\bibitem{ffr:gm}
B. Feigin, E. Frenkel and N. Reshetikhin,
{\it Gaudin model, Bethe ansatz and critical level},
Comm. Math. Phys. {\bf 166} (1994),
27--62.

\bibitem{ffr:oi}
B. Feigin, E. Frenkel and L. Rybnikov,
{\it Opers with irregular singularity and
spectra of the shift of argument
subalgebra},  Duke Math. J. {\bf 155} (2010), 337--363.

\bibitem{fft:gm}
B. Feigin, E. Frenkel and V. Toledano Laredo,
{\it Gaudin models with
irregular singularities},
Adv. Math. {\bf 223} (2010), 873--948.

\bibitem{fjmm:qt}
B. Feigin, M. Jimbo, T. Miwa and E. Mukhin,
{\it Quantum toroidal $\gl_1$-algebra:
plane partitions}, Kyoto J. Math. {\bf 52} (2012), 621--659.

\bibitem{f:lc}
E. Frenkel,
{\it Langlands correspondence for loop groups}, Cambridge
Studies in Advanced Mathematics,
103. Cambridge University Press, Cambridge, 2007.

\bibitem{k:va}
V. Kac,
{\it Vertex algebras for beginners},
University Lecture Series, 10. American Mathematical Society,
Providence, RI, 1997.

\bibitem{m:sf}
{I. G. Macdonald},
{\it Symmetric Functions and Hall Polynomials},
Oxford University Press, Oxford,
1995.

\bibitem{m:ff}
A. I. Molev,
{\it Feigin--Frenkel center in types $B$, $C$ and $D$},
Invent. Math. {\bf 191} (2013), 1--34.

\bibitem{mm:yc}
A. I. Molev and E. E. Mukhin,
{\it Yangian characters and classical $\Wc$-algebras},
in ``Conformal field theory, automorphic forms and related topics"
(W. Kohnen, R. Weissauer, Eds), Springer, 2014, pp. 287--334.

\bibitem{mr:mm}
A. I. Molev and E. Ragoucy,
{\it The MacMahon Master Theorem for right quantum superalgebras
and higher Sugawara operators for $\wh\gl_{m|n}$},
Moscow Math. J. {\bf 14} (2014), 83--119.

\bibitem{fm:qt}
G. S. Mutafyan and B. L. Feigin,
{\it The quantum toroidal algebra $\wh{\wh{\gl_1}}$: calculation of the characters of
some representations as generating functions of plane partitions},
Funct. Anal. Appl. {\bf 47} (2013), 50--61.

\bibitem{oz:sr}
R. C. Orellana and M. Zabrocki,
{\it Some remarks on the characters of the general Lie superalgebra},
{\tt arXiv:math/0008152}.

\bibitem{rt:ip}
M. Ra\"{\i}s and P. Tauvel,
{\it Indice et polyn\^{o}mes invariants pour certaines alg\`{e}bres de Lie},
J. Reine Angew. Math. {\bf 425} (1992), 123--140.

\bibitem{r:si}
L. G. Rybnikov,
{\it The shift of invariants method and the Gaudin model},
Funct. Anal. Appl. {\bf 40} (2006), 188--199.

\bibitem{s:ip}
A. Sergeev,
{\it The invariant polynomials on simple Lie superalgebras},
Represent. Theory {\bf 3} (1999), 250--280.

\bibitem{s:ir}
A. N. Sergeev,
{\it Invariants and representations of classical Lie superalgebras
and their applications to quantum integrable systems} (in Russian),
Doctor of Science Thesis, St. Petersburg Department of
V.A. Steklov Institute of Mathematics of
the Russian Academy of Sciences,
St.-Petersburg, 2008.

\bibitem{sv:dq}
A. N. Sergeev and A. P. Veselov,
{\it Deformed quantum Calogero--Moser problems and Lie superalgebras},
Comm. Math. Phys. {\bf 245} (2004), 249--278.

\bibitem{s:cs}
J. Stembridge,
{\it A characterization of supersymmetric polynomials},
J. Alg. {\bf 95} (1985), 439--444.

\bibitem{sf:fm}
A. V. Stoyanovskii and B. L. Feigin,
{\it Functional models of the representations of current algebras,
and semi-infinite Schubert cells},
Funct. Anal. Appl. {\bf 28} (1994), 55--72.



\end{thebibliography}
\end{document}